\documentclass{amsart}
\usepackage{amsmath,amsfonts,amssymb,amsthm,enumerate}

\newcommand{\R}{{\mathbb R}}
\newcommand{\Z}{{\mathbb Z}}

\newcommand{\C}{{\mathbb C}}

\newcommand{\KK}{{\mathbb K}}

\newcommand{\1}{{\bf 1}}
\newcommand{\F}{{\mathbb F}}

\newcommand{\h}{{\bf h}}

\newcommand{\ch}{{\rm ch }}
\newcommand{\wt}{{\rm wt}}
\newcommand{\Aut}{{\rm Aut}}

\newcommand{\eqa}{\begin{eqnarray}}
\newcommand{\eeqa}{\end{eqnarray}}
\newcommand{\eqn}{\begin{eqnarray*}}
\newcommand{\eeqn}{\end{eqnarray*}}

\newcommand{\Hom}{{\rm Hom}}

\newtheorem{dfn}{Definition}[section]
\newtheorem{pro}[dfn]{Proposition}
\newtheorem{thm}[dfn]{Theorem}

\newtheorem{lem}[dfn]{Lemma}
\newtheorem{cor}[dfn]{Corollary}
\newtheorem{rem}[dfn]{Remark}
\newtheorem{note}[dfn]{Note}

\newcommand{\NO}{\,{\raise0.25em\hbox{$\mathop{\hphantom{\cdot}}\limits^{_{\circ}}_{^{\circ}}$}}\,}

\newcommand{\qe}{\qed\vskip2ex}
\makeatletter
\@addtoreset{equation}{section}

\makeatother

\def\bl{\begin{lem}}
\def\el{\end{lem}}
\def\bt{\begin{thm}}
\def\et{\end{thm}}
\def\bp{\begin{pro}}
\def\ep{\end{pro}}
\def\br{\begin{rem}}
\def\er{\end{rem}}
\def\bc{\begin{cor}}
\def\ec{\end{cor}}
\def\bd{\begin{dfn}\rm}
\def\ed{\end{dfn}}
\def\bn{\begin{note}\rm}
\def\en{\end{note}}
\def\proof{{\it Proof.}}
\def\h{\mathfrak{h}}

\title[On isomorphism problems for VOAs associated with even lattices]{On isomorphism problems for vertex operator algebras associated with even lattices}
\author[H. Shimakura]{Hiroki Shimakura}
\email{shima@auecc.aichi-edu.ac.jp}
\address{Department of Mathematics, Aichi University of Education, 1 Hirosawa, Igaya-cho, Kariya-city, Aichi 448-8542 Japan}

\thanks{The author was partially supported by Grants-in-Aid for Scientific Research (No. 20549004) and Excellent Young Researcher Overseas Visit Program, JSPS}

\begin{document}
\maketitle

\begin{abstract}
In this article, we completely determine the isomorphism classes of lattice vertex operator algebras and the vertex operator subalgebras fixed by a lift of the $-1$-isometry of the lattice. 
We also provide similar results for certain even lattices associated with doubly-even binary codes.
\end{abstract}

\section*{Introduction}
The lattice vertex operator algebra (VOA) $V_L$ associated with an even lattice $L$ is a fundamental example in VOA theory (\cite{Bo,FLM}).
The VOA $V_L$ has an automorphism of order $2$ lifted from the $-1$-isometry of $L$, and the fixed-point subspace $V_L^+$ is a subVOA of $V_L$.
The famous moonshine VOA $V^\natural$ was constructed in \cite{FLM} as an extension of $V_\Lambda^+$ associated with the Leech lattice $\Lambda$.
Therefore, if  $L$ is a sublattice of $\Lambda$, then the VOA $V_L^+$ is contained in $V^\natural$.
Many group theoretical properties of the automorphism group of $V^\natural$, the Monster simple group, can be studied using the subVOA $V_L^+$.
For instance, certain maximal $2$-local subgroups of the Monster simple group were described using $V_L^+$ (\cite{Sh4,Sh}).

It is well-known that the VOAs $V_L$ and $V_N$ are isomorphic if and only if $L$ and $N$ are isomorphic.
However it was shown in \cite{Sh3} that the VOAs $V_{D_{16}^+}^+$ and $V_{E_8^2}^+$ are isomorphic, where $D_{16}^+$ and $E_8^2$ are (non-isomorphic) even unimodular lattices of rank $16$.
Hence it is natural to ask if this case is the only exceptional case.

Given a doubly-even binary code $C$, one can obtain an even lattice $\mathcal{L}(C)=\frac{1}{\sqrt2}\rho^{-1}(C)$ and a certain sublattice $\mathcal{L}^+(C)$ of index $2$, where $\rho$ is the canonical map from $\Z^n$ to $\Z_2^n$.
The Leech lattice $\Lambda$ can be constructed as an overlattice of $\mathcal{L}^+(G_{24})$ (\cite{CS}), where $G_{24}$ is the extended binary Golay code.
The construction of $V^\natural$ in \cite{FLM} is, in some sense, analogous to this construction.
Since there are many analogies among binary codes, lattices and VOAs (\cite{Ho95,Ho,Ho08,Sh}), it is natural to ask if the VOAs $V_L$ and $V_L^+$ would behave like $\mathcal{L}(C)$ and $\mathcal{L}^+(C)$, respectively.

\medskip

In this article, we completely determine when $V_L^+$ and $V_N^+$ are isomorphic as well as when $V_L^+$ and $V_N$ are isomorphic.
We also obtain similar results for the even lattices $\mathcal{L}(C)$ and $\mathcal{L}^+(C)$ associated with doubly-even binary codes $C$.
As an observation, we notice that the isomorphism types of the VOAs $V_L$ and $V_L^+$ are closely related to those of the lattices $\mathcal{L}(C)$ and $\mathcal{L}^+(C)$ (see Table \ref{Table2}).
Note that $\mathcal{C}(K)$ and $\mathcal{C}^+(K)$ in the table are binary codes obtained from Kleinian codes $K$ (\cite{Ho}) and that some results for the even lattices were already given in \cite{KKM}.

\medskip

Next, let us explain our method for the case $V_L^+\cong V_N^+$ (Theorem \ref{TCh2}).
We note that $V_L\cong V_N$ if and only if $L\cong N$ (cf.\ Proposition \ref{T1}), and the case $V_L^+\cong V_N$ (Theorem \ref{T3}) can be solved easily by Theorem \ref{TCh2}.
Clearly, $L\cong N$ implies $V_L^+\cong V_N^+$, and it was shown in \cite{Sh3} that $V_{E_8^2}^+\cong V_{D_{16}^+}^+$.
We assume that $V_L^+\cong V_N^+$ as VOAs.
Let us consider the irreducible simple current module $V_N^-$.
By the assumption, we regard the irreducible $V_N^+$-module $V_N^-$ as an irreducible $V_L^+$-module $M$.
The key is that the properties of $M$ are similar to those of $V_N^-$.
By the classification of the irreducible modules for $V_L^+$ (\cite{DN2,AD}), we have $M\cong V_L^-,\ V_{\lambda+L}^\varepsilon\ (\lambda\notin L,\ \varepsilon\in\{\pm\})$ or $V_L^{T_\chi,\pm}$ as $V_L^+$-modules.

(i): If $M\cong V_L^-$ then, $V_N\cong V_N^+\oplus V_N^-\cong V_L^+\oplus V_L^-\cong V_L$ as VOAs and $L\cong N$.

(ii): Assume $M\cong V_{\lambda+L}^\varepsilon$.
We may assume $\varepsilon=+$ up to conjugation (\cite{Sh2}).
Since the characters of $V_L^-$ and $V_{\lambda+L}^+$ are the same, we obtain an equation about the numbers of norm $2$ vectors in $L$ and $\lambda+L$.
By the characterization of $\mathcal{L}^+(C)$ in \cite{Sh3}, $L\cong\mathcal{L}^+(C)$ and $L+\Z\lambda\cong \mathcal{L}(C)$ for some doubly-even binary code $C$.
Moreover, since $V_{\mathcal{L}(C)}^+\cong V_{\mathcal{L}^+(C)}$ (\cite{FLM}), we obtain $V_L\cong V_{L}^+\oplus V_{L}^-\cong V_L^+\oplus V_{\lambda+L}^+\cong V_N$ as VOAs, and hence $L\cong N$.

(iii): Assume $M\cong V_L^{T_\chi,+}$.
Since the lowest weights are the same, the rank of $L$ is $8$.
Since $M$ is a self-dual simple current, $\sqrt2L^*$ is even.
Hence $L$ contains a sublattice isomorphic to $\sqrt2E_8$.
By the same arguments for $N$, $N$ has the same property.
Since the discriminant groups of $L$ and $N$ are isomorphic, we have $L\cong N$.

(iv): Assume $M\cong V_L^{T_\chi,-}$.
Then the rank of $L$ is $16$.
Moreover, both $L$ and $N$ are unimodular, or $L\cong \mathcal{L}^+(C)$, $N\cong\mathcal{L}^+(D)$ for some doubly-even binary codes $C$ and $D$ with the same weight enumerator.
The first case is an exceptional case.
In the latter case, either both $C$ and $D$ are self-dual, or $C\cong\mathcal{C}^+(K)$, $D\cong\mathcal{C}^+(J)$ for some Kleinian codes $K$ and $J$ with the same weight enumerator.
If both $C$ and $D$ are self-dual then $L\cong \mathcal{L}^+(C)\cong\mathcal{L}^+(D)\cong N$.
By the classification of even Kleinian codes of small length (\cite{Ho}), we have $C\cong \mathcal{C}^+(K)\cong \mathcal{C}^+(J)\cong D$.
Hence $L\cong N$.

This article also contains proofs for the corresponding results in lattices.
These should be useful for understanding the arguments in VOAs since the approach in lattices is quite similar to that in VOAs.

\begin{table}
\caption{Isomorphism problems and the answers}\label{Table2}
\begin{center}
\begin{tabular}{|c| c|}
\hline 
Problem& Answer\\ \hline\hline
$\mathcal{L}(C)\cong \mathcal{L}(D)$ & $C\cong D$ \\ \hline
$V_L\cong V_N$ & $L\cong N$ \\ \hline\hline
$\mathcal{L}^+(C)\cong \mathcal{L}^+(D)$ & $C\cong D$ or $\{C,D\}=\{e_8^2,d_{16}^+\}$ \\ \hline
$V_L^+\cong V_N^+$& $L\cong N$ or $\{L,N\}=\{E_8^2,D_{16}^+\}$\\ \hline
\hline
$\mathcal{L}^+(C)\cong \mathcal{L}(D)$& $C\cong\mathcal{C}(K)$ and $D\cong\mathcal{C}^+(K)$ \\ \hline
$V_L^+\cong V_N$& $L\cong\mathcal{L}(C)$ and $N\cong\mathcal{L}^+(C)$\\ \hline
\end{tabular}
\end{center}
\end{table}

\medskip

\section{Preliminaries}
In this section, we recall or give some definitions and facts required in this article.

\subsection{Kleinian codes}
In this subsection, we recall the basic definitions for Kleinian codes (cf.\ \cite{Ho}).

Denote the elements of the Kleinian four group $\KK\cong \Z_2\times \Z_2$ by $0$, $a$, $b$, $c$, where $0$ is the identity.
A (linear) {\it Kleinian code} $K$ of length $n$ is a subgroup of $\KK^n\cong\Z_{2}^{2n}$.
The dimension of $K$ is $t$ if $K$ has $4^{t}$ elements, where $t\in\frac{1}{2}\Z$.
The {\it weight} $\wt (x)$ of $x=(x_i)\in \KK^n$ is the number of nonzero $x_i$.
For $m\in\Z$, let $K(m)$ denote the set of all elements of weight $m$ in $K$.
A Kleinian code $K$ is {\it even} if $\wt(k)\in2\Z$ for all $k\in K$.
The symmetric bilinear product $\KK^n\times \KK^n\to\F_2$ is defined by $x\cdot y=\sum_{i=1}^nx_i\cdot y_i$, where $a\cdot b=b\cdot a=1$, $a\cdot c=c\cdot a=1$, $c\cdot b=b\cdot c=1$, and zero otherwise.
A Kleinian code is {\it self-dual} if it is equal to its orthogonal complement.
Two Kleinian codes are {\it equivalent} if one of them is obtained from the other by a permutation of the coordinates together with a permutation of the symbols $a,b$ and $c$ at each coordinate.
The {\it weight enumerator} of $K$ is defined by $$W_K(X,Y)=\sum_{k\in K}X^{\wt (k)}Y^{n-\wt (k)}.$$
Up to length $8$, self-dual Kleinian codes were classified in \cite{Ho}.

\bl\label{LHo1} {\rm \cite{Ho}}
\begin{enumerate}
\item There is exactly one even self-dual Kleinian code of length $2$, up to equivalence.
It is equivalent to $\epsilon_2$ generated by $(aa)$ and $(bb)$.
\item There are exactly two even self-dual Kleinian codes of length $4$, up to equivalence.
They are equivalent to $\epsilon_2^2$ and $\delta_4^+$, where $\delta_4^+$ is generated by $(aa00),(a0a0),(a00a)$, and $(bbbb)$
\end{enumerate}
\el

\subsection{Binary codes}
In this subsection, we recall the basic definitions for binary codes (cf.\ \cite{CS}).

A {\it binary ({\it linear}) code} of {\it length} $n$ is a subspace of $\F_2^n$.
The {\it weight} $\wt(x)$ of $x=(x_i)\in\F_2^n$ is the number of nonzero $x_i$.
For a subset $D$ of $\F_2^n$, we denote the set of all elements in $D$ of weight $m$ by $D(m)$.
A binary code $C$ is {\it doubly-even} if $\wt(x)\in4\Z$ for all $x\in C$.
The {\it dual code} $C^\perp$ of $C$ is defined by $C^\perp=\{x\in\F_2^n\mid \langle x,C\rangle=0\}$, where $\langle x,y\rangle=\sum_{i=1}^n x_iy_i$.
A binary code $C$ is {\it self-dual} if $C=C^\perp$.
Two binary codes are {\it equivalent} if one of them is obtained from the other by a permutation of the coordinates.
The {\it weight enumerator} of a coset $x+C\in C^\perp/C$ is defined by $$W_{x+C}(X,Y)=\sum_{y\in x+C}X^{\wt (y)}Y^{n-\wt (y)}=\sum_{m=0}^n |(x+C)(m)|X^mY^{n-m}.$$

We refer to \cite{CS} for the details about the binary codes $e_8$ and $d_{16}^+$.

\bl\label{LCC816} {\rm (cf.\ \cite{CS})}
\begin{enumerate}
\item The extended Hamming code $e_8$ is the unique doubly-even self-dual code of length $8$, up to equivalence.
\item There are exactly two doubly-even self-dual codes of length $16$ , up to equivalence.
They are equivalent to $e_8^2$ and $d_{16}^+$.
\end{enumerate}
\el

\subsection{Lattices}
In this subsection, we recall the basic definitions for lattices (cf.\ \cite{CS}).

Let $(\cdot,\cdot)$ be a positive-definite symmetric bilinear form on $\R^n$.
The {\it norm} of $v\in\R^n$ is $(v,v)$.
For a subset ${U}\subset\R^n$, let us denote by ${U}(m)$ the set of all vectors of norm $m$ in ${U}$.
A subset $L$ of $\R^n$ is a (positive-definite) {\it lattice} of {\it rank} $n$ if there is a basis $e_1,e_2,\dots,e_n$ of $\R^n$ satisfying $L=\bigoplus_{i=1}^n\Z e_i$.
The {\it dual lattice} $L^*$ of $L$ is defined by $L^*=\{u\in \R^n\mid (u,L)\subset\Z\}$.
The {\it discriminant group} of $L$ is the quotient group $L^*/L$.
A lattice $L$ is {\it even} if $( v,v)\in2\Z$ for all $v\in L$, and $L$ is {\it unimodular} if  $L=L^*$.
Two lattices are {\it isomorphic} if one of them is obtained from the other by an orthogonal transformation of $\R^n$.
The {\it theta series} of a coset $\lambda+L\in L^*/L$ is defined by $$\Theta_{\lambda+L}(q)=\sum_{v\in \lambda+L}q^{( v,v)/2}=\sum_{m=0}^\infty |(\lambda+L)(m)|q^{m/2}.$$

We refer to \cite{CS} for the details about the lattices $E_8$ and $D_{16}^+$.

\bl\label{LCL816} {\rm (cf.\ \cite{Wi,CS})}
\begin{enumerate}
\item The root lattice $E_8$ is the unique even unimodular lattice of rank $8$, up to isomorphism.
\item There are exactly two even unimodular lattices of rank $16$, up to isomorphism.
They are isomorphic to $E_8^2$ and $D_{16}^+$.
\end{enumerate}
\el

\subsection{Vertex operator algebras}
In this subsection, we recall the basic definitions for vertex operator algebras (cf.\ \cite{Bo,FLM,FHL}).
Throughout this article, all VOAs are defined over the field $\C$ of complex numbers.

A {\it vertex operator algebra} (VOA) $V$ is a $\Z_{\ge0}$-graded
 vector space $V=\bigoplus_{m\in\Z_{\ge0}}V_m$ equipped with a linear map $Y:V\to ({\rm End}(V))[[z,z^{-1}]]$, $v\mapsto Y(v,z)=\sum_{i\in\Z}a_iz^{-i-1}$ 
and non-zero vectors $\1_V$ and $\omega_V$ satisfying a number of conditions (\cite{Bo,FLM}).
We often denote it by $(V,Y)$ or $V$.
Two VOAs $(V,Y)$ and $(V^\prime,Y')$ are said to be {\it isomorphic} if there exists a linear isomorphism $g$ from $V$ to $V'$ such that $gY(v,z)w=Y'(gv,z)gw$ for all $v,w\in V$ and $g(\omega_V)=\omega_{V'}$.

An (ordinary) module $(M,Y_M)$ for a VOA $V$ is a $\C$-graded vector space $M=\bigoplus_{m\in\C} M_{m}$ equipped with a linear map $Y_M:V\to ({\rm End}(M))[[z,z^{-1}]]$ 
satisfying a number of conditions (\cite{FHL}).
We often denote it by $M$ and its isomorphism class by $[M]$.
The {\it weight} of a homogeneous vector $v\in M_k$ is $k$.
The {\it character} of $M$ is defined by $$\ch_M(q)=q^{-n/24}\sum_{m\in\C}\dim M_{m} q^{m},$$
where $n$ is the central charge of $V$.
A module is {\it self-dual} if its contragredient module (\cite{FHL}) is isomorphic to itself.
If $M$ is irreducible, then there exists $h\in\C$ such that $M=\bigoplus_{m\in\Z_{\ge0}}M_{h+m}$ and $M_h\neq0$.
The number $h$ is called the {\it lowest weight} of $M$.

A VOA $V$ is {\it of CFT type} if $V_0=\C\1$, and $V$ is {\it $C_2$-cofinite} if $V/{\rm Span}_\C\{ a_{-2}b\mid a,b\in V\}$ is finite-dimensional.
A VOA $V$ is {\it simple} if it is an irreducible $V$-module, and $V$ is {\it rational} if any module is completely reducible.
A rational simple VOA $V$ is {\it holomorphic} if any irreducible module is isomorphic to $V$.

The following was proved in \cite{DMb}.
For the lattice VOA $V_L$, see Section 1.5.

\bl\label{LCV816} {\rm \cite{DMb}}
\begin{enumerate}
\item The lattice VOA $V_{E_8}$ is the unique $C_2$-cofinite holomorphic VOA of CFT type of central charge $8$ up to isomorphism.
\item There are exactly two $C_2$-cofinite holomorphic VOAs of CFT type of central charge $16$ up to isomorphism.
They are isomorphic to $V_{E_8^2}$ and $V_{D_{16}^+}$.
\end{enumerate}
\el

Let $M_a,M_b,M_c$ be modules for a simple VOA $V$.
Let $N_{M_a,M_b}^{M_c}$ denote the dimension of the space of all intertwining operators of type $M_a\times M_b\to M_c$, which is called the {\it fusion rule} (\cite{FHL}).
By the definition of the fusion rules, $N_{M_a,M_b}^{M_c}=N_{N_a,N_b}^{N_c}$ if $M_x\cong N_x$ as $V$-modules for all $x=a,b,c$.
Hence, the fusion rules are determined by the isomorphism classes of $V$-modules.
Let $R(V)$ denote the set of all isomorphism classes of irreducible $V$-modules.
For convenience, we use the following expression $$[M_a]\times [M_b]=\sum_{[M]\in R(V)} N_{M_a,M_b}^{M}[M],$$
which is also called the fusion rule.

Let $V(0)$ be a simple VOA.
An irreducible $V(0)$-module $M^1$ is called a {\it simple current} if for any irreducible $V(0)$-module $M^2$, there exists the unique irreducible $V(0)$-module $M^3$ such that $[M^1]\times [M^2]=[M^3]$.
A simple VOA $V$ is called a {\it simple current extension} of $V(0)$ graded by a finite abelian group $A$ if $V$ is the direct sum of non-isomorphic irreducible simple current $V(0)$-modules $\{V(\alpha)\mid \alpha\in A\}$ and for all $\alpha,\beta\in A$, the fusion rule $[V(\alpha)]\times[V(\beta)]=[V(\alpha+\beta)]$ holds.
The uniqueness of simple current extensions was established as follows.

\bp\label{PSC}{\rm \cite[Proposition 4.3]{DM} (cf.\ \cite[Theorem 4.3]{Ho0})} Let $(V,Y)$ be a simple current extension of a simple VOA $V(0)$.
Then the VOA structure of $V$ as a simple current extension of $V(0)$ is uniquely determined by the $V(0)$-module structure of $V$, that is, if $(\tilde{V},\tilde{Y})$ has a VOA structure with $V=\tilde{V}$ as vector spaces and $Y(v,z)=\tilde{Y}(v,z)$ for all $v\in V(0)$, then the VOAs $(V,Y)$ and $(\tilde{V},\tilde{Y})$ are isomorphic.
\ep

When every irreducible $V$-module is a simple current, $\times$ is a binary operation on $R(V)$.
In addition, if $\times$ satisfies the associative law, then $(R(V),\times)$ is a group, which we call the {\it fusion group} of $V$. 

\subsection{Constructions of binary codes, lattices and VOAs}
First, let us consider the binary codes $\mathcal{C}(K)$ and $\mathcal{C}^+(K)$ of length $4m$ obtained from a Kleinian code $K$ of length $m$ \cite[Section 7]{Ho}.

{\it Construction A}: 
\begin{equation}
\mathcal{C}(K)=\hat{K}+d_4^m,\label{Eq:ConstAK}
\end{equation}
where $\ \hat{}\ :K^m\to\F_2^{4m}$ is the map induced from $K\to\F_2^4, 0\mapsto(0000), a\mapsto(1100),b\mapsto (1010),c\mapsto(0110)$, and $d_4^m=\{(0000),(1111)\}^m$.

{\it Construction B}: 
\begin{equation}
\mathcal{C}^+(K)=\hat{K}+(d_4^m)_0,\label{Eq:ConstBK}
\end{equation} 
where $\ \hat{}\ :K^m\to\F_2^{4m}$ is defined as before and $(d_4^m)_0$ is the subcode of $d_4^n$ consisting of vectors of weight divisible by $8$.

\bl\label{LCAB}{\rm (cf.\ \cite[Lemmas 2 and 3]{Ho})} Let $K$ be a Kleinian code of length $m$.
\begin{enumerate}
\item If $K$ is self-dual, then so is $\mathcal{C}(K)$.
\item If $K$ is even, then both $\mathcal{C}(K)$ and $\mathcal{C}^+(K)$ are doubly-even.
\item The weight enumerators of $\mathcal{C}(K)$ and $\mathcal{C}^+(K)$ are given as follows:\begin{eqnarray*}
W_{\mathcal{C}(K)}(X,Y)&=&W_K(X^4+Y^4,2X^2Y^2),\\
W_{\mathcal{C}^+(K)}(X,Y)&=&\frac{1}{2}\left(W_K(X^4+Y^4,2X^2Y^2)+(X^4-Y^4)^m\right).
\end{eqnarray*}
\end{enumerate}
\el

Let us consider some examples.

\bl\label{LC816}
\begin{enumerate}
\item The binary code $\mathcal{C}(\epsilon_2)$ is equivalent to $e_8$.
\item The binary codes $\mathcal{C}(\epsilon_2^2)$ and $\mathcal{C}(\delta_4^+)$ are equivalent to $e_8^2$ and $d_{16}^+$ respectively.
\item The binary codes $\mathcal{C}^+(\epsilon_2^2)$ and $\mathcal{C}^+(\delta_4^+)$ are equivalent.
\end{enumerate}
\el

Next, let us consider the lattices $\mathcal{L}(C)$ and $\mathcal{L}^+(C)$ of rank $n$ obtained from a binary code $C$ of length $n$ (\cite{CS}).
Let $\alpha_1,\alpha_2,\dots,\alpha_n$ be an orthogonal basis of $\R^n$ of norm $2$.

{\it Construction A}: 
\begin{equation}
\mathcal{L}(C)=\sum_{i=1}^n\Z\alpha_i+\sum_{c\in C}\Z\frac{1}{2}\alpha_c,\label{Eq:ConstA}
\end{equation}
where $\alpha_c=\sum_{i=1}^n c_i\alpha_i$ and $c_i\in\{0,1\}$.

{\it Construction B}:
\begin{equation}
\mathcal{L}^+(C)=\sum_{1\le i,j\le n}\Z(\alpha_i+\alpha_j)+\sum_{c\in C}\Z\frac{1}{2}\alpha_c,\label{Eq:ConstB}
\end{equation}
where $\alpha_c$ is defined as before.

\bl\label{LLAB}{\rm \cite{CS}} Let $C$ be a binary code of length $n$.
\begin{enumerate}
\item If $C$ is self-dual, then $\mathcal{L}(C)$ is unimodular.
\item If $C$ is doubly-even, then both $\mathcal{L}(C)$ and $\mathcal{L}^+(C)$ are even.
\item The theta series of $\mathcal{L}(C)$ and $\mathcal{L}^+(C)$ are given as follows:
\begin{eqnarray*}
\Theta_{\mathcal{L}(C)}(q)&=&W_C(\theta_3(q),\theta_2(q)),\\
\Theta_{\mathcal{L}^+(C)}(q)&=&\frac{1}{2}\left(W_C(\theta_3(q),\theta_2(q))+\theta_4(q)^n\right),
\end{eqnarray*}
 where $\theta_3(q)=\sum_{i\in\Z}q^{i^2}$, $\theta_2(q)=\sum_{i\in\Z}q^{(i+1/2)^2}$ and $\theta_4(q)=\sum_{i\in\Z}(-1)^iq^{i^2}$.
\end{enumerate}
\el

Let us consider some examples.
For (3), see \cite[(5.3.1)]{KKM}.

\bl\label{LL816}
\begin{enumerate}
\item The even lattice $\mathcal{L}(e_8)$ is isomorphic to $E_8$.
\item The lattices $\mathcal{L}(e_8^2)$ and $\mathcal{L}(d_{16}^+)$ are isomorphic to $E_8^2$ and $D_{16}^+$ respectively.
\item The lattices $\mathcal{L}^+(e_8^2)$ and $\mathcal{L}^+(d_{16}^+)$ are isomorphic.
\end{enumerate}
\el

Finally, we consider the VOAs $V_L$ and $V_L^+$ of central charge $n$ obtained from an even lattice $L$ of rank $n$ (\cite{Bo,FLM}).

Set $\mathfrak{h}_L = \C \otimes_{\Z} L$ and extend $(\cdot,\cdot)$ to a symmetric $\C$-bilinear form on $\h_L$.
We view $\mathfrak{h}_L = \C \otimes_{\Z} L$ as an abelian Lie algebra.
Let $\hat{\mathfrak{h}}_L = \mathfrak{h}_L \otimes \C [t,t^{-1}] \oplus \C {\bf c}$ be its affine Lie algebra. 
Set ${\hat{\frak{h}}_L}^-=\frak{h}_L\otimes t^{-1}\C[t^{-1}]$ and let $S(\hat{\frak{h}}^-_L)$ be the symmetric algebra of $\hat{\frak{h}}^-_L$.
Then $S(\hat{\frak{h}}^-_L)\cong\C[\alpha \otimes t^n\mid\alpha \in \mathfrak{h}_L, n < 0]\cdot \1$ is the unique irreducible $\hat{\mathfrak{h}}_L$-module such that ${\bf c}$ acts as $1$ and $\alpha \otimes t^n\cdot \1 = 0$ if $\alpha \in \mathfrak{h}_L$ and $n \ge 0$.

Let us consider the twisted group algebra of the additive group $L$.
Let $\langle\kappa_L\rangle$ be a cyclic group of order $2$ and let
\begin{eqnarray}
1\to \langle\kappa_L\rangle  \to \hat{L}\ \overset{\rho}{\to} L\ \to 1\label{Seq:rho}
\end{eqnarray}
be a central extension of $L$ by $\langle\kappa_L\rangle$ with the commutator map $c_0(\alpha, \beta) = (\alpha, \beta) \pmod 2$, $\alpha, \beta \in L$.
Then we obtain the twisted group algebra $\C\{L\}= \C[\hat{L}]/(\kappa_L+1)$, where $\C[\hat{L}]$ is the usual group algebra of the group $\hat{L}$.
The {\it lattice VOA} $V_L$ associated with $L$ is defined to be $V_L=S(\hat{\frak{h}}^-_L) \otimes \C\{L\}$ as a vector space.
We refer to \cite{FLM} for the definition of the vertex operator on $V_L$.
Note that the central charge of $V_L$ is $n$.

Let $\theta_{V_L}$ be an automorphism of $V_L$ of order $2$ which is a lift of the $-1$-isometry of $L$.
Let $V_L^+$ denote the subspace of $V_L$ fixed by $\theta_{V_L}$.
Then $V_L^+$ is a vertex operator subalgebra of $V_L$.
Note that the automorphism $\theta_{V_L}$ is not unique but the isomorphism type of $V_L^+$ is independent of the choice of $\theta_{V_L}$ \cite[Appendix D]{DGH}.

As a summary, we obtain two VOAs of central charge $n$:
$$V_L=S(\hat{\frak{h}}^-_L) \otimes \C\{L\}\quad {\rm and}\quad V_L^+=\{v\in V_L\mid \theta_{V_L}(v)=v\}.$$
We call the constructions of $V_L$ and $V_L^+$ from an even lattice $L$ {\it Construction A} and {\it Construction B} respectively.

\bl\label{LVAB} {\rm \cite{Do,FLM}} Let $L$ be an even lattice of rank $n$.
\begin{enumerate}
\item If $L$ is unimodular, then $V_L$ is holomorphic.
\item The characters of $V_L$ and $V_L^+$ are given as follows:
\begin{eqnarray*}
\ch_{V_L}(q)&=&\frac{\Theta_L(q)}{\eta(q)^n},\\
\ch_{V_L^+}(q)&=&\frac{1}{2}\left(\frac{\Theta_L(q)}{\eta(q)^n}+\frac{\eta(q)^n}{\eta(2q)^n}\right),
\end{eqnarray*}
where $\eta(q)=q^{1/24}\Pi_{i=1}^\infty (1-q^i)$.
\end{enumerate}
\el

The following lemma was verified in {\rm \cite[Lemma 3.4]{Sh3}.

\bl\label{LV816} The VOAs $V_{E_8^2}^+$ and $V_{D_{16}^+}^+$ are isomorphic.
\el

Now, we recall some facts about irreducible modules for $V_L^+$ (\cite{FLM,DN2,AD}).
Let $V_L^-$ be the $-1$-eigenspace of $\theta_{V_L}$ in $V_L$.
Then $V_L^-$ is an irreducible $V_L^+$-module.
For a coset $\lambda+L\in L^*/L$, $V_{\lambda+L}=S(\hat{\frak{h}}^-_L) \otimes \C[\lambda+L]$ is an irreducible $V_L$-module.
If $\lambda\notin L/2$, then $V_{\lambda+L}$ is also irreducible for $V_L^+$.
If $\lambda\in L/2$, then $V_{\lambda+L}$ is the direct sum of two irreducible $V_L^+$-submodules $V_{\lambda+L}^+$ and $V_{\lambda+L}^-$.
Set $J=\{a^{-1}\theta_{V_L}(a)\mid a\in \hat{L}\}$.
Then $J$ is a normal subgroup of $\hat{L}$.
Let $T_\chi$ be an irreducible module for $\hat{L}/J$ with character $\chi$ satisfying $\chi(\kappa J)=-1$.
Then there are two irreducible $V_L^+$-modules $V_L^{T_\chi,\pm}$ associated with $T_\chi$.

It was shown in \cite{DN2,AD} that any irreducible $V_L^+$-module is isomorphic to one of $V_{\lambda+L}^\pm$ $(\lambda\in L^*\cap (L/2))$, $V_{\mu+L}$ $(\mu\in L^*\setminus (L/2))$ and $V_L^{T_{\chi},\pm}$.
By the fusion rules for irreducible $V_L^+$-modules \cite{Ab,ADL}, one can obtain the following lemma.

\bl\label{LADL} The irreducible $V_L^+$-module $V_{\mu+L}$ $(\mu\in L^*\setminus (L/2))$ is not a simple current.
\el

The characters of $V_{\lambda+L}^\pm$ and $V_L^{T_{\chi},\pm}$ are given as follows (\cite{FLM}):

\begin{eqnarray}
\ch_{V_{\lambda+L}^\pm}(q)&=&\frac{1}{2}\Big(\frac{\Theta_{\lambda+L}(q)}{\eta(q)^n}\Big)\quad (\lambda\notin L),\label{Eq:untwist}\\
\ch_{V_L^{T_{\chi},\pm}}(q)&=&\frac{\dim T_\chi}{2}\Big(\frac{\eta(q)^n}{\eta(q^{1/2})^n}\pm\frac{\eta(q^2)^n\eta(q^{1/2})^n}{\eta(q)^{2n}}\Big).\label{Eq:twist}
\end{eqnarray}

\section{Isomorphism problems for the even lattices $\mathcal{L}(C)$ and $\mathcal{L}^+(C)$}

\subsection{A characterization of doubly-even binary codes obtained from Kleinian codes}
In this subsection, we characterize the doubly-even binary codes obtained by Construction B from Kleinian codes, which will play an important role in Sections 2.3 and 2.4.
Note that this is a code analogue of the characterization of the even lattices obtained by Construction B (see Proposition \ref{PSh}) given in \cite{Sh3}.

Let $u_i=(0^{4(i-1)}1^40^{n-4i})\in\F_2^n$ for $1\le i\le n/4$.
If $m=n/4\in\Z$, then $d_4^m={\rm Span}_{\F_2}\{u_i\mid 1\le i\le m\}$ and $(d_4^m)_0={\rm Span}_{\F_2}\{u_i+u_j\mid 1\le i,j\le m\}$.
For our purpose, we need the following lemma.

\bl\label{LCC} Let $C$ be a doubly-even binary code of length $n=4m\in4\Z$.
\begin{enumerate}
\item If $d_4^m\subset C$, then $C\cong\mathcal{C}(K)$ for some even Kleinian code $K$ of length $m$.
\item Assume that $(d_4^m)_0\subset C$, $u_1\notin C$ and $u_1\in C^\perp$.
Then $C\cong\mathcal{C}^+(K)$ and $C+\F_2 u_1\cong\mathcal{C}(K)$ for some even Kleinian code $K$ of length $m$.
\end{enumerate}
\el
\proof\ (1): Let $x=(x_1,x_2,\dots, x_m)\in C$, where $x_i\in\F_2^4$.
Since $C$ is doubly even, $\langle x,u_i\rangle=\langle x_i,(1111)\rangle=0$, equivalently, $\wt(x_i)\in2\Z$.
It follows from $\hat{\KK}+\F_2(1111)=\{v\in\F_2^4\mid \wt(v)\in2\Z\}$ that $x_i\in\hat{\KK}$ or $x_i+(1111)\in\hat{\KK}$.
For the map $\hat{}$\ , see (\ref{Eq:ConstAK}).
Hence there is a unique $k_x\in {\KK}^{m}$ such that $x\in\hat{k}_x+d_4^{m}$.
Set $K=\{k_x\mid x\in C\}$.
Then $C\subset\mathcal{C}(K)$.
It follows from $d_4^m\subset C$ that $\hat{k}_x\in C$.
Thus by (\ref{Eq:ConstAK}) $C\cong\mathcal{C}(K)$.
Since $C$ is a doubly-even binary code, $K$ is an even Kleinian code.

(2): Set $\tilde{C}=C+\F_2u_1$.
It follows from $u_1\in C^\perp$ and $\wt(u_1)=4$ that $\tilde{C}$ is doubly-even.
Since $u_1\notin C$, we have $C\subsetneq\tilde{C}$.
Let $y=(y_1,y_2,\dots,y_m)\in C^\perp\setminus\tilde{C}^\perp$, where $y_i\in\F_2^4$.
It follows from $\langle y,u_i\rangle\neq 0$ and $\langle y,u_i+u_j\rangle=0$ for all $i,j$ that $\wt(y_i)=1$ or $3$ for all $i$.
Since $\tilde{C}$ is doubly-even, $\tilde{C}\subset\tilde{C}^\perp$.
Hence $u_i\in\tilde{C}^\perp$ for any $i$, and we may assume that $\wt(y_i)=1$ for all $i$, and $y=(00010001\dots0001)$ up to coordinate permutation.

Let $x=(x_1,x_2,\dots,x_m)\in C$, where $x_i\in\F_2^4$.
Since $u_i\in C^\perp$ for all $i$, we have $\langle x_i,(1111)\rangle=0$.
Hence there is $k_x\in\KK^m$ such that $x\in \hat{k}_x+d_4^m$.
Since $\langle x,y\rangle=\sum_{i=1}^m\langle x_i,(0001)\rangle=0$ and $\langle \hat{k}_x,y\rangle=0$, we have $x\in \hat{k}_x+(d_4^m)_0$.
By the same argument as in (1), $K=\{k_x\mid x\in C\}$ is an even Kleinian code and $C\cong \mathcal{C}^+(K)$ (cf.\ (\ref{Eq:ConstBK})).
Moreover, $C+\F_2 u_1\cong\mathcal{C}(K)$.
\qe

\bp\label{MTC} Let $C$ be a doubly-even binary code of length $n$.
If there exists a coset $x+C\in C^\perp/C$ satisfying $|(x+C)(4)|\ge n/4+|C(4)|$ then $C\cong\mathcal{C}^+(K)$ and $C+\F_2x\cong\mathcal{C}(K)$ for some even Kleinian code $K$ of length $n/4$.
\ep
\proof\ In order to apply Lemma \ref{LCC} (2), we will show that the coset $x+C$ contains $\{u_i\mid 1\le i\le n/4\}$.
Assume $\{u_i\mid 1\le i\le r\}\subset x+C$ and $r<n/4$.
We claim that the coset $x+C$ contains $(0^{4r},z)$ with $\wt(z)=4$.

Let $v=(v_1,v_2,\dots,v_{r},z)\in (x+C)(4)$, where $v_i\in\F_2^4$ and $z\in\F_2^{n-4r}$.
Then $v+u_i\in C$.
It follows from $u_i\in C^\perp$ that $\langle v+u_i,u_i\rangle=\langle v_i,(1111)\rangle=0$, equivalently, $\wt(v_i)\in2\Z$.
If $\wt(v_i)=4$ for some $i$, then $v=u_i$ since $\wt(v)=4$.
Set $$Y=\{(y_1,y_2,\dots, y_r,z)\in (x+C)(4)\mid \wt(y_i)=2\ {\rm for\ some}\ 1\le i\le r\}.$$
For $y\in Y$, let $m(y)=\min\{i\mid \wt(y_i)=2\}$.
In order to give an upper bound for $|Y|$, we consider the map $\psi:Y\to C(4)$, $y\mapsto y+u_{m(y)}$ and show its injectivity.
Assume $\psi(y)=\psi(y')$.
Then $y+u_{m(y)}=y'+u_{m(y')}$, equivalently, $y+y'=u_{m(y)}+u_{m(y')}$.
By the definition of $u_i$, $\wt(u_{m(y)}+u_{m(y')})=0$ or $8$.
If $\wt(u_{m(y)}+u_{m(y')})=8$, then $m(y)\neq m(y')$.
However, by $\wt(y+y')=8$, we have $\wt(y_{m(y)})=\wt(y_{m(y')})=\wt(y'_{m(y)})=\wt(y'_{m(y')})=2$, and $m(y)=m(y')$, which a contradiction.
Hence $\wt(u_{m(y)}+u_{m(y')})=0$, equivalently, $y=y'$.
Thus $\psi$ is injective.
In particular, $|Y|\le |C(4)|$.

By the definition of $Y$, we have
\begin{equation*}\{(0,\dots,0,z)\in (x+C)(4)\mid z\in\F_2^{n-4r}\}=(x+C)(4)\setminus ({
Y}\cup\{u_i\mid 1\le i\le r\}).\label{Eq:x+C}\end{equation*}
By the upper bound of $|Y|$, $$|(x+C)(4)\setminus ({
Y}\cup\{u_i\mid 1\le i\le r\})|\ge |(x+C)(4)|-|C(4)|-r,$$
and it is greater than $0$ by $r<n/4$ and the claim follows.
Hence there is an element $(0^{4r},z)\in x+C$ with $\wt(z)=4$.
We define it to be $u_{r+1}$ up to coordinate permutation.
Thus the coset $x+C$ contains $\{u_i\mid 1\le i\le n/4\}$ up to coordinate permutation, and $n/4$ must be an integer.
Moreover, $u_i+u_j\in C$ for all $i,j$.
Hence $u_1\in C^\perp$, $u_1\notin C$, and ${\rm Span}_{\F_2}\{u_i+u_j\mid 1\le i,j\le m\}=(d_4^m)_0\subset C$.
Therefore the proposition follows from Lemma \ref{LCC} (2).
\qe

\br If $C= \mathcal{C}^+(K)$, then $|((1^40^{n-4})+C)(4)|= n/4+|C(4)|$.
\er

\subsection{Even lattices $\mathcal{L}(C)$}
In this subsection, we study the isomorphism classes of $\mathcal{L}(C)$.
The following proposition is well-known (cf.\ \cite[Chapter 18, Proposition 4]{CS}).
However, we give a sketch of a proof since it would be useful for understanding the arguments in the VOA case (see Section 3.1).

\bp\label{PLA} Let $C$ and $D$ be doubly-even binary codes of length $n$.
Then $\mathcal{L}(C)\cong \mathcal{L}(D)$ if and only if ${C}\cong{D}$.
\ep
\proof\ The ``if" part is obvious.
For the ``only" part, we assume that $\mathcal{L}(C)\cong \mathcal{L}(D)$.
One can verify that the Weyl group of the norm $2$ vectors in $\mathcal{L}(C)$ is transitive on the set of all orthogonal bases of norm $2$ in $\mathcal{L}(C)$.
Hence we may assume that the canonical bases $\{\alpha_i\mid 1\le i\le n\}$ of $\mathcal{L}(C)$ and $\mathcal{L}(D)$ in (\ref{Eq:ConstA}) are the same up to automorphism.
Since the binary code $C$ can be recovered from $\{\alpha_i\mid 1\le i\le n\}$, we obtain $C\cong D$.\qed

\subsection{Even lattices $\mathcal{L}^+(C)$}
In this subsection, we study the isomorphism classes of $\mathcal{L}^+(C)$.

The following lemma can be easily proved by using the quadratic spaces of plus type over $\F_2$ associated to doubly-even binary codes of length $8k$ containing $(1^{8k})$.

\bl\label{LC8k} For any doubly-even binary code of length $8k$, there exists a doubly-even self-dual binary code of length $8k$ containing it.
\el

By (\ref{Eq:ConstBK}), $(d_4^m)_0\subset \mathcal{C}^+(K)$ and $(00010001\dots0001)\in \mathcal{C}^+(K)^\perp$.
Hence by \cite[Lemma 3.3.2]{KKM}, we obtain the following lemma.

\bl\label{LKKM332} Let $K$ be an even Kleinian code of length $m$.
Let $\rho$ be the orthogonal transformation of $\R^{4m}$ defined by
\begin{eqnarray*}
\rho(\alpha_{j})=\left\{\begin{array}{cc}
 \mbox{$\frac{1}{2}(\alpha_{4i-3}+\alpha_{4i-2}+\alpha_{4i-1}+\alpha_{4i})$} & \mbox{${\rm if}\ j=4i$},\\
 \mbox{$\frac{1}{2}(-\alpha_{4i-3}-\alpha_{4i-2}-\alpha_{4i-1}+\alpha_{4i})+\alpha_j$} & \mbox{${\rm if}\ 4i-3\le j\le 4i-1$},\end{array}
\right.
\end{eqnarray*}
where $\{\alpha_i\mid 1\le i\le 4m\}$ is the basis of $\R^{4m}$ as defined in (\ref{Eq:ConstA}) and (\ref{Eq:ConstB}).
Then $\rho$ is an isomorphism of lattices from $\mathcal{L}(\mathcal{C}^+(K))$ to $\mathcal{L}^+(\mathcal{C}(K))$.
\el

The following proposition was shown in \cite{KKM} if $n>16$, and the case $n\le 16$ is a finite problem.
However, we give a proof here since it is useful for understanding the VOA case (Theorem \ref{TCh2}).

\bt\label{MTT} {\rm (cf.\ \cite[Theorem 3]{KKM})} Let $C$ and $D$ be doubly-even binary codes of length $n$.
Then $\mathcal{L}^+(C)\cong\mathcal{L}^+(D)$ if and only if one of the following holds.
\begin{enumerate}
\item $C$ and $D$ are equivalent.
\item $C$ and $D$ are doubly-even self-dual binary codes of length $16$.
\end{enumerate}
\et
\proof\ The ``if" part follows from Lemmas \ref{LCC816} (2) and \ref{LL816} (3).
For the ``only" part, we assume that $\mathcal{L}^+(C)\cong \mathcal{L}^+(D)$.
Then by Lemma \ref{LLAB} (3) and the fact that $\theta_2$ and $\theta_3$ are algebraically independent (\cite[Section 2.9]{Eb}), $W_C(X,Y)=W_D(X,Y)$.
Since $\mathcal{L}^+(C)\cong \mathcal{L}^+(D)$, there is an isomorphism $\varphi:\mathcal{L}^+(D)^*/\mathcal{L}^+(D)\to \mathcal{L}^+(C)^*/\mathcal{L}^+(C)$ of their discriminant groups.
Let $\lambda+\mathcal{L}^+(C)=\varphi(\alpha_1+\mathcal{L}^+(D))$, where $\{\alpha_i\mid 1\le i\le n\}$ is the canonical basis of $\R^n$ as defined in (\ref{Eq:ConstA}) and (\ref{Eq:ConstB}).
Without loss of generality, we assume $\langle \lambda,\lambda\rangle=2$.
By (\ref{Eq:ConstB}), $$\mathcal{L}^+(C)^*=\mathcal{L}(C^\perp)+\Z\frac{1}{4}\alpha_{(1^n)}.$$
For $c=(c_i)\in \F_2^n$, let $\varepsilon_c$ be the orthogonal transformation defined by $\varepsilon_c(\alpha_i)=(-1)^{c_i}\alpha_i$.
Then $\varepsilon_c$ is an automorphism of $\mathcal{L}^+(C)$ if $c\in C^\perp$.
By the actions of $\varepsilon_c$, we may assume that $\lambda=\alpha_1$, $\alpha_y/2$ $(y\in C^\perp(4)\setminus C)$, $\alpha_{(1^n)}/4-\alpha_1$ or $\alpha_{(1^n)}/4$.

(i) Assume that $\lambda=\alpha_1$.
Then $\mathcal{L}(C)\cong\mathcal{L}(D)$, and by Proposition \ref{PLA}, $C\cong D$.

(ii) Assume that $\lambda=\alpha_y/2$.
Comparing the numbers of vectors of norm $2$ in $\alpha_1+\mathcal{L}^+(D)$ and $\lambda+\mathcal{L}^+(C)$, we obtain $2^3|D(4)|+2n=2^3|(y+C)(4)|$.
By $W_C(X,Y)=W_D(X,Y)$, we have $|C(4)|=|D(4)|$.
Combining the equations, we obtain $$2^3|C(4)|+2n=2^3|(y+C)(4)|.$$
Then by Proposition \ref{MTC}, $C\cong \mathcal{C}^+(K)$ and $C+\F_2y\cong\mathcal{C}(K)$ for some even Kleinian code $K$ of length $n/4$.
Hence
\begin{eqnarray*}
\mathcal{L}(D)=\mathcal{L}^+(D)\cup (\alpha_1+\mathcal{L}^+(D))\cong \mathcal{L}^+(C)\cup (\lambda+\mathcal{L}^+(C))\cong\mathcal{L}^+(\mathcal{C}(K))\cong\mathcal{L}(C),
\end{eqnarray*}
where the last isomorphism is given in Lemma \ref{LKKM332}.
Thus by Proposition \ref{PLA}, $C\cong D$.

(iii) Assume that $\lambda=\alpha_{(1^n)}/4-\alpha_1$.
Since the norm of $\lambda$ is $2$, $n$ must be $8$.
Since $2\lambda\in \mathcal{L}^+(C)$, $C$ contains $(1^8)$.
It follows from $W_C(X,Y)=W_D(X,Y)$ that $D$ also has $(1^8)$ and $\dim C=\dim D$.
One easily see that there are only four doubly-even binary codes of length $8$ containing $(1^8)$ up to equivalence, and those are determined by the dimensions.
Hence $C\cong D$.

(iv) Assume that $\lambda=\alpha_{(1^n)}/4$.
Since the norm of $\lambda$ is $2$, $n$ must be $16$.
Let $k$ be the dimension of $C$.
Then $C^\perp/C\cong \Z_2^{16-2k}.$
By $|(\alpha_1+\mathcal{L}^+(C))(2)|=|(\lambda+\mathcal{L}^+(C))(2)|$, we obtain $2^3|C(4)|+2^5=2^k$, equivalently, \begin{equation}|C(4)|=2^{k-3}-4.\label{Eq:C4}\end{equation}
If $C$ is self-dual, then so is $D$, and hence (2) holds.
Assume that $C$ is not self-dual.
Then by $|C(4)|\ge0$, we have $5\le k\le 7$.
By Lemma \ref{LC8k}, there is a doubly-even self-dual binary code $E$ of length $16$ such that $C\subsetneq E$.
By Lemma \ref{LCL816} (2), $E\cong d_{16}^+$ or $e_8^2$, and $|E(4)|=28$.
Hence there exists a coset $x+C\in E/C$ such that
$$|(x+C)(4)|\ge\frac{|E(4)|-|C(4)|}{|E/C|-1}=\frac{2^5-2^{k-3}}{2^{8-k}-1}=2^{k-3}=|C(4)|+\frac{16}{4}.$$
By Proposition \ref{MTC}, $C\cong\mathcal{C}^+(K)$ for some even Kleinian code $K$ of length $4$.
By $|\mathcal{C}^+(K)|=2^3|K|$, we have $|K|=2^{k-3}$.
By Lemma \ref{LCAB} (3) and (\ref{Eq:C4}), $|K(2)|=|C(4)|/2=2^{k-4}-2$.
Since $|K(4)|=|K|-|K(2)|-|K(0)|=2^{k-4}+1$, we obtain
\begin{equation}
W_K(X,Y)=(2^{k-4}+1)X^4+(2^{k-4}-2)X^2Y^2+Y^4.\label{Eq:WK}
\end{equation}

Replacing $C$ by $D$, we also have $C\cong D$, or $D\cong\mathcal{C}^+(J)$ for some even Kleinian code $J$ of length $4$ with the same properties as $K$.
One can verify that there are only four even Kleinian codes of length $4$ satisfying (\ref{Eq:WK}) up to equivalence (cf.\ \cite{Ho}).
Two of them are characterized by the weight enumerators, and the others are the two even self-dual Kleinian codes $\epsilon_2^2$ and $\delta_4^+$ of length $4$ (Lemma \ref{LHo1} (2)).
By Lemma \ref{LC816} (3) $\mathcal{C}^+(\epsilon_2^2)\cong\mathcal{C}^+(\delta_4^+)$.
Thus we obtain $C\cong D$.\qe

\subsection{Even lattices $\mathcal{L}(C)$ and $\mathcal{L}^+(C)$}
In this subsection, we study the isomorphism classes of  $\mathcal{L}(C)$ and $\mathcal{L}^+(C)$.

\bt\label{TAB} Let $C$ and $D$ be doubly-even binary codes of length $n$.
Then $\mathcal{L}^+(C)\cong \mathcal{L}(D)$ if and only if there exists an even Kleinian code $K$ such that $C\cong \mathcal{C}(K)$ and $D\cong \mathcal{C}^+(K)$.
\et
\proof\ The ``if" part follows from Lemma \ref{LKKM332}.
For the ``only" part, we assume that $\mathcal{L}^+(C)\cong \mathcal{L}(D)$.
Comparing the numbers of vectors of norm $2$ in the lattices (Lemma \ref{LLAB} (3)), we obtain
$$8|C(4)|=2n+16|D(4)|.$$
Since $\mathcal{L}^+(C)\cong \mathcal{L}(D)$, there is an isomorphism $\varphi:\mathcal{L}^+({C})^*/\mathcal{L}^+(C)\to \mathcal{L}({D})^*/\mathcal{L}({D})$ of their discriminant groups.
By $\mathcal{L}(D)^*=\mathcal{L}(D^\perp)$, there is $x\in D^\perp\setminus D$ such that $\alpha_x/2+\mathcal{L}(D)=\varphi(\alpha_1+\mathcal{L}^+(C))$.
Comparing the numbers of vectors of norm $2$ in $\alpha_1+\mathcal{L}^+(C)$ and $\alpha_x/2+\mathcal{L}(D)$, we obtain
$$2n+8|C(4)|=16|(x+D)(4)|.$$
Combining the equations, we obtain 
$$|(x+D)(4)|=n/4+|D(4)|.$$
By Proposition \ref{MTC}, $D\cong \mathcal{C}^+(K)$ for some even Kleinian code $K$.
Then by Proposition \ref{LKKM332} and an assumption,  
$$\mathcal{L}^+(\mathcal{C}(K))\cong\mathcal{L}(\mathcal{C}^+(K))\cong \mathcal{L}(D)  \cong \mathcal{L}^+(C).$$
By Theorem \ref{MTT}, $C\cong \mathcal{C}(K)$, or both $C$ and $\mathcal{C}(K)$ are (non-isomorphic) doubly-even self-dual binary codes of length $16$.
In the first case, we are done.
In the latter case,  Lemmas \ref{LCC816} (2) and \ref{LC816} (2) show $C\cong \mathcal{C}(K')$, where $\{K,K'\}=\{\epsilon_2^2,\delta_4^+\}$.
By Lemma \ref{LC816} (3), $\mathcal{C}^+(K)\cong\mathcal{C}^+(K')$.
Hence we have $D\cong \mathcal{C}^+(K')$ and $C\cong \mathcal{C}(K')$.
\qe

\section{Isomorphism problems for the VOAs $V_L$ and $V_L^+$}

\subsection{VOAs $V_L$}
In this subsection, we study the isomorphism classes of $V_L$.
The following proposition is well-known.
Its proof is similar to that of Proposition \ref{PLA}.

\bp\label{T1} Let $L$ and $N$ be even lattices of rank $n$.
Then $V_L\cong V_N$ if and only if $L\cong N$.
\ep
\proof\ The ``if" part is obvious.
For the ``only" part, we assume $V_L\cong V_N$.
By the VOA structure of $V_L$ (\cite{FLM}), $S(\hat{\frak{h}}^-_L)$ is a subVOA and $V_L\cong \bigoplus_{\alpha\in L}S(\hat{\frak{h}}^-_L)\otimes e^\alpha$ as $S(\hat{\frak{h}}^-_L)$-modules.
Note that $\frak{h}_L(-1)\1=S(\hat{\frak{h}}^-_L)\cap (V_L)_1$.
Then for $h\in \frak{h}_L$, $S(\hat{\frak{h}}^-_L)\otimes e^\alpha$ is the eigenspace of the $0$-th product of $h(-1)\1$ with eigenvalue $(\alpha,h)$.
By the same argument, $V_N\cong\bigoplus_{\beta\in N}S(\hat{\frak{h}}^-_N)\otimes e^\beta$ as $S(\hat{\frak{h}}^-_N)$-modules, and $\frak{h}_N(-1)\1=S(\hat{\frak{h}}^-_N)\cap (V_N)_1$.
Since $\dim (V_L)_0=1$, $(V_L)_1$ has a finite dimensional reductive Lie algebra structure by the $0$-th product.
Using the fact that $V_L\cong V_N$, we may view both $\frak{h}_L(-1)\1$ and $\frak{h}_N(-1)\1$ as Cartan subalgebras of $(V_L)_1$.
Hence there exists $g\in\langle\exp(v_0)\mid v\in (V_L)_1\rangle\subset\Aut(V_L)$ such that $g(\frak{h}_N(-1)\1)=\frak{h}_L(-1)\1$.
Since $S(\hat{\frak{h}}^-_N)$ is generated by $\frak{h}_N(-1)\1$ as a VOA, we have $g(S(\hat{\frak{h}}^-_N))=S(\hat{\frak{h}}^-_L)$.
Hence for any $\beta\in N$, $g(S(\hat{\frak{h}}^-_N)\otimes e^\beta)=S(\hat{\frak{h}}^-_L)\otimes e^\alpha$ for some $\alpha\in L$.
Set $\bar{g}(\beta)=\alpha$.
Since $g$ is an automorphism of a VOA, $\bar{g}$ is a linear isomorphism from $N$ to $L$, and $(\bar{g}(\beta),\bar{g}(\beta))=(\beta,\beta)$ for all $\beta\in N$.
Thus $\bar{g}$ is an orthogonal transformation, and $N\cong L$.\qe

\subsection{VOAs $V_L^+$} 
In this subsection, we study the isomorphism classes of $V_L^+$.

The following proposition, which is similar to Proposition \ref{MTC}, was proved in \cite{Sh3}.

\bp\label{PSh}{\rm \cite[Theorem 2.2]{Sh3}} Let $L$ be an even lattice of rank $n$.
If there exists $\lambda\in (L^*\cap L/2)$ such that $|(\lambda+L)(2)|\ge2n+|L(2)|$ then $L\cong\mathcal{L}^+(C)$ and $L+\Z\lambda\cong\mathcal{L}(C)$ for some doubly-even binary code $C$.
\ep

Next, we will consider an even lattice $L$ such that $\sqrt2L^*$ is even.

\bl\label{L2ele} Let $L$ be an even lattice.
If $\sqrt2 L^*$ is even, then $2L^*\subset L$.
\el 
\proof\ Since $\sqrt2L^*$ is even, we have $\sqrt2L^*\subset (\sqrt2L^*)^*=L/\sqrt2$, and $2L^*\subset L$.\qe

\br In \cite{Sh2,Sh3}, an even lattice $L$ was said to be $2$-elementary totally even if $2L^*\subset L$ and if $\sqrt2 L^*$ is even.
By the lemma above, if $\sqrt2 L^*$ is even, then $L$ is $2$-elementary totally even.
\er

The following lemma can be proved by using the quadratic spaces over $\F_2$ associated to $2$-elementary totally even lattices of rank $8k$.
Note that it is of plus type by \cite[Theorem 1]{CS2}.

\bl\label{LL8k} Let $L$ be an even lattice of rank $8k$ such that $\sqrt2L^*$ is even.
Then there exists an even unimodular lattice of rank $8k$ containing $L$. 
\el

Let us review some properties of irreducible modules for $V_L^+$.

\bl\label{L1} The following three conditions are equivalent.
\begin{enumerate}
\item Any irreducible $V_L^+$-module is a self-dual simple current.
\item An irreducible $V_L^+$-module $V_L^{T_\chi,\varepsilon}$ is self-dual.
\item $\sqrt2L^*$ is even.
\end{enumerate}
\el
\proof\ Obviously, (1) shows (2).
By Lemma \ref{L2ele} and \cite{ADL}, one can see that (3) implies (1).

Assume (2) and let $V_L^{T_\chi,\varepsilon}$ be a self-dual irreducible $V_L^+$-module.
Let $M$ be the contragredient module of $V_L^{T_\chi,\varepsilon}$.
Then by \cite[Proposition 3.7]{ADL}, $M\cong V_L^{T_{\chi^\prime},\varepsilon}$ and $\chi^\prime(s)=(-1)^{( \rho{(s)},\rho{(s)})/2}\chi(s)$ for $s\in Z({\hat{L}/J})$, where $\rho$ is given in (\ref{Seq:rho}).
Since $V_L^{T_\chi,\varepsilon}$ is self-dual, we have $\chi^\prime=\chi$.
It follows from $\rho({Z(\hat{L}/J}))=2L^*/2L$ that $( v,v)\in 4\Z$ for all $v\in 2L^*$.
Hence $\sqrt2L^*$ is even, and we obtain (3).\qe

Note that the fusion rules of $V_L^+$ are associative \cite[Theorem 5.18]{ADL}.

\bl\label{LFg} {\rm (cf.\ \cite{ADL})} Let $L$ be an even lattice such that $\sqrt2L^*$ is even.
Then the fusion group of $V_L^+$ is isomorphic to an elementary abelian $2$-group of order $2^{k+2}$, where $k$ is the integer satisfying $L^*/L\cong\Z_2^k$.
\el

A triality automorphism of $V_{\mathcal{L}^+(C)}^+$ which sends $V_{\mathcal{L}^+(C)}^-$ to $V_{\alpha_1+\mathcal{L}^+(C)}^+$ was given in \cite[Chapter 11]{FLM}. 
The following lemma can be verified directly.

\bl\label{LFLM} {\rm (cf.\ \cite[Chapter 11]{FLM})} Let $C$ be a doubly-even binary code.
Then there is an isomorphism of VOAs from $V_{\mathcal{L}^+(C)}$ to $V_{\mathcal{L}(C)}^+$.
\el

Now, we prove the following.
Its proof is similar to that of Theorem \ref{MTT}.

\bt\label{TCh2} Let $L$ and $N$ be even lattices of rank $n$.
Then the VOAs $V_L^+$ and $V_N^+$ are isomorphic if and only if one of the following holds.
\begin{enumerate}
\item $L$ and $N$ are isomorphic.
\item $L$ and $N$ are even unimodular lattices of rank $16$.
\end{enumerate}
\et
\proof\ The ``if" part follows from Lemmas \ref{LCL816} (2) and \ref{LV816}.
For the ``only" part, we assume that $V_L^+\cong V_N^+$.
By Lemma \ref{LVAB} (2), $\Theta_L(q)=\Theta_N(q)$, and hence $\ch V_L^-=\ch V_N^-$.
We regard the irreducible $V_N^+$-module $V_N^-$ as an irreducible $V_L^+$-module.
Since $V_N^-$ is a self-dual simple current, it is isomorphic to $V_L^-$, $V_{\lambda+L}^\varepsilon$ ($\lambda\notin L$) or $V_L^{T_\chi,\varepsilon}$ as $V_L^+$-modules under the isomorphism of $V_N^+$ and $V_L^+$ by Lemma \ref{LADL}.

(i) Assume that $V_N^-\cong V_L^-$.
Then $V_N\cong V_L^+\oplus V_L^-$ is a simple current extension of $V_L^+$.
Hence by Proposition \ref{PSC} $V_N\cong V_L$ as VOAs, and by Proposition \ref{T1} $L\cong N$.

(ii) Assume that $V_N^-\cong V_{\lambda+L}^\varepsilon$.
Since $\Hom(L^*/2L^*,\Z_2)\subset\Aut(V_L^+)$ (\cite{FLM}), we may assume that $\varepsilon=+$ up to conjugation (\cite{Sh2}).
It follows from $\ch V_L^-=\ch V_N^-$ that $\ch V_L^-=\ch V_{\lambda+L}^+$.
Comparing the dimensions of the weight $1$ subspaces of $V_L^-$ and $V_{\lambda+L}^+$ (Lemma \ref{LVAB} and (\ref{Eq:untwist})), we obtain $|L(2)|+2n=|(\lambda+L)(2)|$.
Hence by Proposition \ref{PSh} $L\cong\mathcal{L}^+(C)$ and $L+\Z\lambda\cong\mathcal{L}(C)$ for some doubly-even binary code $C$.
By Proposition \ref{PSC}, as VOAs, 
$$V_N\cong V_L^+\oplus V_N^-\cong V_{L}^+\oplus V_{\lambda+L}^+\cong V_{\mathcal{L}(C)}^+\cong V_L,$$
where the last isomorphism is given in Lemma \ref{LFLM}.
By Proposition \ref{T1}, we obtain $L\cong N$.

(iii) Assume that $V_N^-\cong V_L^{T_\chi,-}$.
By Lemma \ref{LVAB} and (\ref{Eq:twist}), the lowest weights of $V_N^-$ and $V_L^{T_\chi,-}$ are $1$ and $n/16+1/2$, respectively.
Hence $n=8$.
Since $V_N^-$ is a self-dual, so is $V_L^{T_\chi,-}$.
Hence $\sqrt2L^*$ is even by Lemma \ref{L1}.
Since $\sqrt2(\sqrt2L^*)^*=L$ is even, we have $\sqrt2L^*\subset E_8$ by Lemmas \ref{LCL816} (1) and \ref{LL8k}.
Hence $E_8\subset (\sqrt2L^*)^*=L/\sqrt2$, equivalently, $\sqrt2E_8\subset L$.
Replacing $N$ by $L$, we have (i), (ii) or (iii) for $V_L^-$.
The first two cases show $L\cong N$.
Hence we may assume that $\sqrt2E_8\subset N$ by the same argument.
One easily checks that there are only five even overlattices of $\sqrt2E_8$ up to isomorphism, and those are determined by the discriminant groups.
Since the fusion groups of $V_L^+$ and $V_N^+$ are isomorphic, so are the discriminant groups of $L$ and $N$ by Lemma \ref{LFg}.
Thus $L\cong N$.

(iv) Assume that $V_N^-\cong V_L^{T_\chi,+}$.
By (\ref{Eq:twist}), the lowest weight of $V_L^{T_\chi,+}$ is $n/16$.
By the same argument as in (iii), $\sqrt2L^*$ is even and $n=16$.
Comparing the dimensions of the weight $1$ subspaces of $V_N^-$ and $V_L^{T_\chi,+}$ (Lemma \ref{LVAB} and (\ref{Eq:twist})), we obtain 
\begin{equation}
|L(2)|=2^{9-k}-32,\label{Eq:L2}
\end{equation}
where $|L^*/L|=2^{2k}$.
If $L$ is unimodular, then so is $N$, and hence (2) holds.
Assume that $L$ is not unimodular.
By Lemma \ref{LL8k} there exists an even unimodular lattice $P$ of rank $16$ such that $L\subsetneq P$.
By Lemma \ref{LCL816} (2), $|P(2)|=480$.
Hence there exists a coset $\mu+L\in P/L$ such that $$|(\mu+L)(2)|\ge \frac{|P(2)|-|L(2)|}{|P/L|-1}=\frac{512-2^{9-k}}{2^k-1}=2^{9-k}=|L(2)|+2\times 16.$$
By Proposition \ref{PSh}, $L\cong \mathcal{L}^+(C)$ for some $t$-dimensional doubly-even binary code $C$ of length $16$.
It follows from $|\mathcal{L}^+(C)^*/\mathcal{L}^+(C)|=2^{18-2t}$ that $t=9-k$.
By Lemma \ref{LLAB} and $n=16$, $|\mathcal{L}^+(C)(2)|=2^3|C(4)|+2\times 16$.
Hence, by (\ref{Eq:L2}), we obtain 
\begin{equation*}
|C(4)|=2^{t-3}-4.
\end{equation*}
Since $|C(4)|\ge0$, we have $5\le t\le 8$.
By the same arguments as in (iv) of Theorem \ref{MTT} (see (\ref{Eq:C4})), if $C$ is not self-dual, then $C\cong\mathcal{C}^+(K)$ for some even Kleinian code $K$ of length $4$ with the weight enumerator (\ref{Eq:WK}).

Replacing $L$ by $N$, we also obtain $L\cong N$ or $N\cong\mathcal{L}^+(D)$ for some doubly-even binary code $D$ of length $16$ with the same properties as $C$.
If $C$ is self-dual then so is $D$, and $L\cong N$ by Lemma \ref{LL816} (3).
If not, $C\cong\mathcal{C}^+(K)$ and $D\cong\mathcal{C}^+(J)$ for some even Kleinian codes $K$ and $J$ of length $4$. 
By Lemma \ref{LFg}, we have $L^*/L\cong N^*/N$, $\dim C=\dim D$, and $\dim K=\dim J$.
Then by (\ref{Eq:WK}), both Kleinian codes $K$ and $J$ have the same weight enumerator.
Thus $K\cong J$ or $\{K,J\}=\{\varepsilon_2^2,\delta_4^+\}$ (cf.\ \cite{Ho}).
Therefore by Lemma \ref{LC816} $C\cong D$, and $L\cong \mathcal{L}^+(C)\cong\mathcal{L}^+(D)\cong N$.\qe

\subsection{VOAs $V_L$ and $V_L^+$}
In this subsection, we study the isomorphism classes of $V_L$ and $V_L^+$.
The proof of the following theorem is similar to that of Theorem \ref{TAB}.

\bt\label{T3} Let $L$ and $N$ be even lattices of rank $n$.
Then the VOAs $V_L^+$ and $V_N$ are isomorphic if and only if
there exists a doubly-even binary code $C$ such that $L\cong\mathcal{L}(C)$ and $N\cong\mathcal{L}^+(C)$.
\et
\proof\ Lemma \ref{LFLM} proves the ''if" part.
For the ``only" part, we assume that $V_L^+\cong V_N$.
Comparing the dimensions of the weight $1$ subspaces of $V_L^+$ and $V_N$ (Lemma \ref{LVAB} (2)), we obtain 
$$\frac{|L(2)|}{2}=n+|N(2)|.$$
Let $V_{\lambda+N}$ be an irreducible $V_N$-module such that $V_{\lambda+N}\cong V_L^-$ under the isomorphism of $V_N$ and $V_L^+$.
Comparing the dimensions of the weight $1$ subspaces of $V_L^-$ and $V_{\lambda+N}$, we obtain 
$$\frac{|L(2)|}{2}+n=|(\lambda+N)(2)|.$$
Combining the equations, we obtain
$$|N(2)|+2n=|(\lambda+N)(2)|.$$
In addition, the fusion rule $[V_L^-]\times[V_L^-]=[V_L^+]$ (\cite{ADL}) shows  $[V_{\lambda+N}]\times[V_{\lambda+N}]=[{V_N}]$.
Hence $2\lambda\in N$.
By Proposition \ref{PSh}, $N\cong \mathcal{L}^+(C)$ for some doubly-even binary code $C$ of length $n$.
By Lemma \ref{LFLM}, we obtain $V_L^+\cong V_N\cong V_{\mathcal{L}^+(C)}\cong V_{\mathcal{L}(C)}^+$ as VOAs.
By Theorem \ref{TCh2}, $L\cong \mathcal{L}(C)$ or both $L$ and $\mathcal{L}(C)$ are (non-isomorphic) even unimodular lattices of rank $16$.
In the first case, we are done.
In the latter case, Lemmas \ref{LCL816} (2) and \ref{LL816} (2) show $L\cong \mathcal{L}(D)$, where $\{C,D\}=\{e_8^2,d_{16}^+\}$.
By Theorem \ref{MTT}, $\mathcal{L}^+(C)\cong \mathcal{L}^+(D)$.
Hence $N\cong \mathcal{L}^+(D)$ and $L\cong \mathcal{L}(D)$.
\qe

\paragraph{\bf Acknowledgements.} The author thanks the referee for helpful suggestions and useful comments.
He also thanks Professor Ching Hung Lam for helpful comments.
Part of the work was done when he was visiting Imperial College London in 2010.
He thanks the staff for their help.

\end{document}